\documentclass[10pt,leqno]{article}

\usepackage{amsmath,amssymb,amsthm,mathrsfs,dsfont,mathtools}

\usepackage[margin=3cm,headheight=14pt]{geometry}

\usepackage{titlesec}
\usepackage{color}
\usepackage{microtype}
\usepackage{indentfirst}
\usepackage{enumitem}
\usepackage{needspace}
\usepackage{cite}
\usepackage[hidelinks]{hyperref}
\usepackage{fancyhdr}
\titleformat{\subsection}{\it}{\thesubsection.\enspace}{1pt}{}

\newtheorem{theorem}{Theorem}[section]
\newtheorem{lemma}[theorem]{Lemma}
\newtheorem{definition}[theorem]{Definition}
\newtheorem{corollary}[theorem]{Corollary}
\newtheorem{proposition}[theorem]{Proposition}
\newtheorem{remark}[theorem]{Remark}
\numberwithin{equation}{section}

\allowdisplaybreaks
\setlist[enumerate,1]{label=(\roman*),ref=\roman*,leftmargin=2.6em,itemsep=4pt,topsep=5pt,parsep=0pt}

\newcommand{\R}{\mathbb{R}}
\newcommand{\dd}{\,\mathrm{d}x}
\newcommand{\ind}{\mathbf{1}}
\newcommand{\uh}{u_h}
\newcommand{\RR}{\mathcal{R}}
\newcommand{\norm}[2]{\left\|#1\right\|_{#2}}
\newcommand{\step}[1]{\par\medskip\Needspace{4\baselineskip}%
  \noindent\textbf{#1}\par\nobreak\smallskip\noindent\ignorespaces}
\DeclareMathOperator{\Div}{div}
\DeclareMathOperator{\curl}{curl}
\DeclareMathOperator{\logp}{log_+}

\hypersetup{
  pdftitle={Liouville theorems for the stationary Navier-Stokes equations under partial component conditions},
  pdfauthor={Cuiyang Liu and Wei Luo},
  pdfsubject={Stationary Navier-Stokes equations; fixed velocity projections; Lorentz spaces},
  pdfkeywords={Liouville theorem, partial component conditions, two components, one component, Lorentz spaces}}

\begin{document}
\title{Liouville theorems for the stationary Navier--Stokes equations under partial component conditions.
\hspace{-4mm}}
\author{Cuiyang $\mbox{Liu}^1$\footnote{E-mail: liucy25@mail2.sysu.edu.cn},\quad
Wei $\mbox{Luo}^1$\footnote{E-mail: luow23@mail.sysu.edu.cn}\\
$^1\mbox{Department}$ of Mathematics,
Sun Yat-sen University, Guangzhou 510275, China}
\date{}
\maketitle
\hrule

\begin{abstract}
We prove Liouville theorems for smooth stationary incompressible
Navier--Stokes flows in $\R^3$ with finite Dirichlet integral and
velocity vanishing uniformly at infinity. If
the partial components of velocity $(u_1,u_2)$ belong to $L^{4,\infty}$, then we show that $v$ must vanish. The second result only requires that
one component $u_1$ belongs to $L^{3/2,\infty}$. The proofs are based on the head-pressure sign
and a localized energy identity, combined with different pressure
estimates for the two cases.
\end{abstract}
\par\medskip

\vspace*{10pt}
\tableofcontents

\section{Introduction}
\subsection{The stationary Liouville problem}
We study the stationary incompressible Navier--Stokes equations
\begin{equation}\label{eq:NS}
    -\Delta u+(u\cdot\nabla)u+\nabla p=0,
    \qquad \Div u=0
    \quad\text{in }\R^3,
\end{equation}
where $u=(u_1,u_2,u_3)$ is the velocity and $p$ is the pressure.
The viscosity is normalized to one, and there is no external force.
Throughout the paper, the background assumptions are
\begin{equation}\label{eq:background}
    \int_{\R^3}|\nabla u|^2\dd<\infty,
    \qquad
    \lim_{R\to\infty}\sup_{|x|>R}|u(x)|=0.
\end{equation}

The study of steady viscous flow has long been tied to finite-Dirichlet-integral
solutions and their behavior at infinity. Leray's work
\cite{Leray1933} laid the foundations for the existence theory of
stationary flows. The subsequent analysis of exterior problems made
clear the importance of distinguishing a finite Dirichlet integral
from a prescribed velocity at infinity. In two dimensions, this issue
was developed in particular by Gilbarg and Weinberger
\cite{GilbargWeinberger1978}. A systematic account of the steady-state
theory, including existence, uniqueness, and asymptotic properties, is
given in Galdi's monograph \cite{Galdi2011}.
In the three-dimensional whole space, a central question is whether
\eqref{eq:background} alone forces a smooth solution of \eqref{eq:NS}
to vanish; see \cite{Galdi2011,CoiculescuYang2026}.
The main obstruction is an energy flux at infinity. The background
assumptions yield $u\in L^6$, but this does not by itself justify
testing the equations globally with $u$: both the convective flux
and the pressure flux have to be controlled.

One line of progress strengthens the integrability or decay of the
full velocity. The classical condition $u\in L^{9/2}(\R^3)$ closes
the localized energy estimate; see \cite{Galdi2011}.
Chae and Wolf \cite{ChaeWolf2016} weakened this condition by a
logarithmic factor. Seregin \cite{Seregin2016} used the Stokes system
with a divergence-free drift to obtain criteria involving Morrey-space
control. Seregin and Wang \cite{SereginWang2019} formulated sufficient
conditions through Lorentz norms on large annuli, while Chamorro,
Jarr\'in and Lemari\'e-Rieusset \cite{ChamorroJarrinLemarie2021}
developed further criteria in Lebesgue and Morrey spaces. These works
refine the same basic question: what information at large distances
is sufficient to make the localized energy flux disappear? In weak
Lorentz spaces, finiteness of a global quasi-norm does not generally
imply that its restrictions to exterior regions tend to zero
\cite{Hunt1966,BennettSharpley1988}. This endpoint issue is relevant to the conditions considered here.

A second line of work uses quantities that retain more of the
structure of the equations. Seregin \cite{Seregin2018} obtained
Liouville theorems from mean-oscillation bounds on a vector potential
of the velocity. Coiculescu and Yang \cite{CoiculescuYang2026}
developed conditional results based on antiderivatives and line
integrals of the velocity. The head pressure, namely the sum of the
pressure and the kinetic energy density, provides another useful
quantity. Its elliptic equation and sign have been used in
\cite{Chae2014,ChaeWolf2016}; more recent conditions compare the
behavior of the head pressure with that of the velocity or the
gradient of the head pressure \cite{Chae2026}. These approaches show that a
Liouville criterion need not be stated solely as an integrability
assumption on $|u|$.

This leads to a more specific question: how much of the velocity
must be controlled? Phan \cite{Phan2020} studied weighted mixed-norm
conditions. Chae \cite{Chae2023} and Zhang and Zu \cite{ZhangZu2025}
considered anisotropic conditions on velocity components. Related
results of Vergara-Hermosilla \cite{Vergara2026} impose conditions
on the spherical radial component $u(x)\cdot x/|x|$. These settings
must be distinguished: mixed norms separate spatial variables,
whereas a component assumption selects part of the vector field;
moreover, a radial direction varies with $x$ and cannot be made into
a fixed Cartesian direction by a constant rotation. In the present
paper, we keep \eqref{eq:background} and impose an ordinary,
unweighted weak-Lorentz condition on either two fixed components
or one fixed component. The two cases require different pressure
estimates but share the same localized energy framework.

\subsection{One-component regularity criteria for the time-dependent Navier--Stokes equations}
For the time-dependent Navier--Stokes equations, conditions on one
velocity component have mainly been used to prevent a loss of
regularity. For a fixed unit vector $e$, a representative condition is
\begin{equation}\label{eq:evolution-component}
    \int_0^T
       \norm{u(t)\cdot e}{\dot H^{1/2+2/\gamma}(\R^3)}^\gamma
       \,\mathrm{d}t<\infty,
\end{equation}
where $\dot H^s$ denotes the homogeneous $L^2$ Sobolev space.
Under their respective initial-data and solution assumptions,
Chemin and Zhang \cite{CheminZhang2016} established such a criterion
for $4<\gamma<6$, Chemin, Zhang and Zhang
\cite{CheminZhangZhang2017} extended the range to $4<\gamma<\infty$,
and Han, Lei, Li and Zhao \cite{HanLeiLiZhao2019} treated the full
range $2\le\gamma<\infty$.
The distinction from our results is the conclusion as well as the
hypothesis. Those results give regularity or continuation from a
space-time Sobolev condition; here the solution is already smooth,
and a spatial Lorentz condition forces it to vanish. In particular,
applying a regularity criterion to a time-independent solution does
not by itself prove a Liouville theorem.

\subsection{Two main theorems}
Our purpose is to find conditions on only part of the velocity that
force $u=0$.
We write that
\begin{equation}\label{eq:coordinates}
    x=(x_h,z),\quad x_h=(x_1,x_2),\quad z=x_3,
    \quad \uh=(u_1,u_2),\quad w=u_3.
\end{equation}

\begin{theorem}[Two-component condition]\label{thm:two}
Let $(u,p)$ be a smooth solution of \eqref{eq:NS} satisfying
\eqref{eq:background}. If
\begin{equation}\label{eq:two-condition}
    \uh\in L^{4,\infty}(\R^3),
\end{equation}
then $u=0$ and $p$ is constant.
\end{theorem}

The second theorem imposes the integrability condition on only one
velocity component.

\begin{theorem}[One-component condition]\label{thm:one}
Let $(u,p)$ be a smooth solution of \eqref{eq:NS} satisfying
\eqref{eq:background}. If
\begin{equation}\label{eq:one-condition}
    u_3\in L^{3/2,\infty}(\R^3),
\end{equation}
then $u=0$ and $p$ is constant.
\end{theorem}

By orthogonal invariance, $\uh$ can be replaced by the pair of
components along any two fixed orthonormal directions, and $u_3$
by the component along any fixed unit vector. The coordinate
transformation is recorded at the beginning of
Section~\ref{sec:common}. Thus the use of $\uh$ and $u_3$ does not
restrict the orientation of either condition.

The terms one-component and two-component refer only to the additional
hypotheses \eqref{eq:one-condition} and \eqref{eq:two-condition}.
Both background conditions in \eqref{eq:background} concern the full
velocity. No further integrability is required of the uncontrolled
components. Neither theorem assumes axial symmetry or independence
of a spatial variable. The two criteria balance the number of
controlled components against their integrability; neither is asserted
to contain every previously known anisotropic condition.

\begin{remark}\label{rem:comparison}
\leavevmode\par
\begin{enumerate}
\item Since $u$ is bounded under \eqref{eq:background}, the
$L^{3/2,\infty}$ condition is strictly stronger than the
$L^{4,\infty}$ condition on the same component. Thus controlling
fewer components comes at the cost of stronger integrability.

\item The one-component proof uses the stronger integrability to close
a scalar momentum estimate and obtain $p\in L^1$. The
$L^{4,\infty}$ bound does not close the same estimate, so the
two-component proof instead exploits incompressibility to control
the pressure flux directly.

\item No smallness condition is imposed on the norm of $u$ in the main results.
\end{enumerate}
\end{remark}

Section~\ref{sec:consequences} records the consequences for smaller
exponents and conditions at infinity.

\subsection{Main ideas of the proofs}
We describe the common structure of the proofs. The background
assumptions give $u\in L^6\cap L^\infty$. After subtracting a unique
constant, we continue to denote the pressure by $p$; it satisfies
\begin{equation}\label{eq:intro-canonical}
    p=\sum_{i,j=1}^3\RR_i\RR_j(u_i u_j)\in L^3,
\end{equation}
where $\RR_j$ denotes the $j$th Riesz transform. The head pressure
$Q=p+|u|^2/2$ is nonpositive. Therefore
\begin{equation}\label{eq:intro-sign}
    q:=-p\ge\frac12|u|^2,
    \qquad |Q|\le q.
\end{equation}
The pressure limit at infinity follows from Galdi's theorem, and
the sign is the standard maximum-principle consequence of the
head-pressure equation; see Lemmas~\ref{lem:normalization}
and~\ref{lem:sign}. It allows pressure estimates to control the
kinetic and convective terms in the same local energy identity.

For the two-component condition, incompressibility gives
\begin{equation}\label{eq:intro-cancellation}
    -\Delta p
       =\sum_{\alpha,\beta=1}^2
           \partial_\alpha\partial_\beta(u_\alpha u_\beta)
          +2\partial_z(\uh\cdot\nabla_hw).
\end{equation}
The purely vertical quadratic term has disappeared. Differentiating
the associated potential representation gives
$\nabla p\in L^{4/3,2}$, and hence
\begin{equation}\label{eq:intro-two-estimates}
    p\in L^{12/5,2},\qquad p/\lvert x_h\rvert\in L^{4/3,2}.
\end{equation}
The second conclusion is a cylindrical Hardy estimate. The remaining
difficulty is that $L^{4/3,2}$ does not pair directly with
$L^{4,\infty}$ in $L^1$. A local $L^{4,2}$ estimate has a square-root
logarithmic loss. Square-summability of the weighted pressure on
disjoint cylindrical annuli removes this loss along a subsequence.
The same energy identity, now on cylinders with $H=R^{3/2}$, completes
the proof.

For the one-component condition, the third momentum equation is tested
against compactly supported approximations to
\begin{equation*}
    \psi_\ell(x)=\frac{\ell z}{\sqrt{\ell^2+|x|^2}}.
\end{equation*}
The derivative $\partial_z\psi_\ell$ is positive, the horizontal
convection can be absorbed by the pressure term, and
$\Delta\psi_\ell$ is uniformly bounded in $L^{3,1}$.
Cylindrical cutoffs with radius $R$ and height scale $H=\delta R^4$
justify the test through the ordered limits
$R\to\infty$, $\delta\downarrow0$, and $\ell\to\infty$.
The resulting estimate is
\begin{equation}\label{eq:intro-one-estimate}
    \norm{p}{L^1}
       \le C\left(\norm{w}{L^{3/2,\infty}}+\norm{w}{L^2}^2\right).
\end{equation}
Here $w\in L^2$ follows by interpolation with $L^6$. The pressure
sign then gives $u\in L^2$, and the energy identity yields $u=0$.

Thus the one-component argument obtains finite total pressure by a
scalar momentum test, whereas the two-component argument controls the
pressure flux directly. Both treat the full weak endpoint without
assuming that the weak norm of the controlled projection tends to zero
at infinity. The exponents $3/2$ and $4$ are endpoints of the respective
estimates, not claims of optimality for the stationary equations.

Section~\ref{sec:common} establishes the common analytic framework.
Section~\ref{sec:two} proves Theorem~\ref{thm:two}, and
Section~\ref{sec:one} proves Theorem~\ref{thm:one}. Section~\ref{sec:consequences} collects exterior
and pointwise consequences and explains the endpoint restrictions.

\section{Preliminaries}\label{sec:common}
Unless a domain is specified, all integrals and function spaces are
taken over $\R^3$. We use the coordinates in \eqref{eq:coordinates},
write $\nabla_h=(\partial_1,\partial_2)$, and set $\partial_z=\partial_3$.
Constants denoted by $C$ may change from line to line. They are independent
of all cut off scales; dependence on other parameters will be indicated.


\subsection{Lorentz spaces and standard estimates}
We collect the standard facts used below. The estimates in this
subsection are stated without proof; see
\cite{Hunt1966,ONeil1963,Stein1970,BennettSharpley1988}.

\begin{definition}[Lorentz spaces]\label{def:lorentz}
Let $f$ be a measurable scalar or vector field on $\R^3$. Its
distribution function and decreasing rearrangement are, respectively,
\begin{equation*}
    d_f(\lambda)=|\{x\in\R^3:|f(x)|>\lambda\}|,
    \qquad
    f^*(t)=\inf\{\lambda\ge0:d_f(\lambda)\le t\},
\end{equation*}
for $\lambda\ge0$ and $t>0$, with $\inf\varnothing=\infty$.
Thus $f^*$ is the nonincreasing rearrangement of $|f|$ on $(0,\infty)$.
For $1\le r<\infty$ and $1\le s<\infty$, define
\begin{equation}\label{eq:lorentz-definition}
    \begin{aligned}
    \norm{f}{L^{r,s}}
       &=\left(\int_0^\infty[t^{1/r}f^*(t)]^s
                         \frac{\mathrm{d}t}{t}\right)^{1/s},\\
    \norm{f}{L^{r,\infty}}
       &=\sup_{t>0}t^{1/r}f^*(t)
        =\sup_{\lambda>0}\lambda d_f(\lambda)^{1/r}.
    \end{aligned}
\end{equation}
The Lorentz space $L^{r,s}(\R^3)$ consists of the measurable fields
for which the corresponding quasi-norm is finite. Restriction to a
measurable set $E$ means multiplication by $\ind_E$.
\end{definition}

\begin{lemma}[H\"older's estimates
{\cite{Hunt1966,ONeil1963,BennettSharpley1988}}]\label{lem:lorentz-holder}
Let $1<a,b<\infty$, $1\le r<\infty$, and
$1\le s,s_1,s_2\le\infty$ satisfy
$1/r=1/a+1/b$ and $1/s=1/s_1+1/s_2$, with $1/\infty=0$.
For measurable $f$ and $g$ with finite right-hand side,
\begin{equation}\label{eq:product}
    \norm{fg}{L^{r,s}}
       \le C\norm{f}{L^{a,s_1}}\norm{g}{L^{b,s_2}}.
\end{equation}
The constant depends only on the exponents. The case $r=s=1$
gives the integral H\"older inequality, since $L^{1,1}=L^1$.
\end{lemma}
The following lemma is used to deal with the pressure term in Lorentz spaces.
\begin{lemma}[Potential estimates
{\cite{ONeil1963,Stein1970,Hunt1966}}]\label{lem:lorentz-operators}
Let $I_1=(-\Delta)^{-1/2}$ on $\R^3$. For $1<r<3$,
$1\le s\le\infty$, and $1/r^*=1/r-1/3$,
\begin{equation}\label{eq:potential}
    \norm{I_1f}{L^{r^*,s}}\le C\norm{f}{L^{r,s}}.
\end{equation}
For $1<r<\infty$, $1\le s\le\infty$, and $1\le i,j\le3$,
\begin{equation*}
    \norm{\RR_jf}{L^{r,s}}
       +\norm{\RR_i\RR_jf}{L^{r,s}}
       \le C\norm{f}{L^{r,s}}.
\end{equation*}
Here $\RR_j$ is the Riesz transform with Fourier multiplier
$-\mathrm{i}\xi_j/|\xi|$, and the constants depend only on the exponents.
\end{lemma}
In order to obtain the localization of solutions we need the restrictions and tails properties in Lorentz spaces.
\begin{lemma}[
{\cite{Hunt1966,BennettSharpley1988}}]\label{lem:lorentz-basic}
Let $1<r<\infty$.
\begin{enumerate}
\item For $1\le s\le\infty$, replacing $f^*$ in
\eqref{eq:lorentz-definition} by
\begin{equation*}
    f^{**}(t)=\frac1t\int_0^t f^*(\tau)\,\mathrm{d}\tau
\end{equation*}
defines a Banach norm $\norm{f}{L^{r,s}}^\#$ satisfying
\begin{equation*}
    \norm{f}{L^{r,s}}
       \le\norm{f}{L^{r,s}}^\#
       \le\frac{r}{r-1}\norm{f}{L^{r,s}}.
\end{equation*}
\item One has $L^{r,r}=L^r$, with equality of the quasi-norm in
\eqref{eq:lorentz-definition} and the usual $L^r$ norm. If
$1\le s_1\le s_2\le\infty$, then
\begin{equation*}
    \norm{f}{L^{r,s_2}}\le C\norm{f}{L^{r,s_1}}.
\end{equation*}
In particular, $L^{r,s}\subset L^{r,\infty}$.
\item For $f\in L^{r,\infty}$, $0<a<r$, and a measurable set $E$
of finite measure,
\begin{equation}\label{eq:finite-volume}
    \int_E|f|^a\dd
       \le\frac{r}{r-a}\norm{f}{L^{r,\infty}}^a
                |E|^{1-a/r}.
\end{equation}
\item If $1\le s<\infty$, $f\in L^{r,s}$, and
$\ind_{E_n}\to0$ almost everywhere, then
$\norm{f\ind_{E_n}}{L^{r,s}}\to0$. In particular,
\begin{equation*}
    \lim_{R\to\infty}
       \norm{f\ind_{\{|x|>R\}}}{L^{r,s}}=0.
\end{equation*}
This tail property is not valid for a general function in
$L^{r,\infty}$.
\end{enumerate}
\end{lemma}

When Minkowski's inequality is used, we use the equivalent Banach
norm in Lemma~\ref{lem:lorentz-basic}(i).
Distribution-function calculations always use the quasi-norm in
Definition~\ref{def:lorentz}.

\subsection{Pressure normalization and its sign}
We recall Galdi's pressure limit at infinity
\cite[Theorem~X.5.1, p.~688]{Galdi2011}, together with the
standard whole-space pressure representation
\cite[Remark~1.1]{Chae2014}. The estimate below is the
Calder\'on--Zygmund estimate in Lemma~\ref{lem:lorentz-operators}.

\begin{lemma}\label{lem:normalization}
Let $(u,p)$ be a smooth solution of \eqref{eq:NS} satisfying
\eqref{eq:background}. Then $u\in L^6\cap L^\infty$. There is a
unique constant $\bar p$ such that $p-\bar p\in L^3$, and
\begin{equation}\label{eq:canonical}
    p-\bar p=\sum_{i,j=1}^3\RR_i\RR_j(u_i u_j),
    \qquad \norm{p-\bar p}{L^3}\le C\norm{u}{L^6}^2.
\end{equation}
Moreover, $p(x)\to\bar p$ uniformly as $|x|\to\infty$.
\end{lemma}

From now on, we subtract the constant in \eqref{eq:canonical} and
continue to write $p$ for the normalized pressure. Thus
$p=\sum_{i,j=1}^3\RR_i\RR_j(u_i u_j)\in L^3$ and $p(x)\to0$
as $|x|\to\infty$. Taking the divergence of \eqref{eq:NS} gives
\begin{equation}\label{eq:poisson}
    -\Delta p=\sum_{i,j=1}^3\partial_i\partial_j(u_i u_j)
       =\sum_{i,j=1}^3\partial_i u_j\partial_j u_i.
\end{equation}
We set
\begin{equation}\label{eq:Qq}
    Q=p+\frac12|u|^2,\qquad q=-p.
\end{equation}
Then $Q\in L^3$ and $Q(x)\to0$ uniformly as $|x|\to\infty$.
No integrability of $p$ below $L^3$ is assumed at this stage.
The following standard head-pressure identity and sign property
are recalled from \cite[Section~2, (2.1)--(2.4)]{Chae2026};
see also \cite[Section~2.2]{ChaeWolf2016}.

\begin{lemma}\label{lem:sign}
Under the assumptions of Lemma~\ref{lem:normalization},
\begin{equation}\label{eq:bernoulli}
    -\Delta Q+u\cdot\nabla Q=-|\curl u|^2,
\end{equation}
and
\begin{equation}\label{eq:sign}
    Q\le0,\qquad q\ge\frac12|u|^2\ge0,
    \qquad |Q|\le q=|p|.
\end{equation}
\end{lemma}

\subsection{Energy localization}
\begin{lemma}\label{lem:energy}
Under the assumptions of Lemma~\ref{lem:normalization}, every
$\phi\in C_c^\infty(\R^3)$ satisfies
\begin{equation}\label{eq:local-energy}
    \int|\nabla u|^2\phi\dd
       =\frac12\int|u|^2\Delta\phi\dd
          +\int Q\,u\cdot\nabla\phi\dd,
\end{equation}
where all integrals are absolutely convergent.
If, in addition, the normalized pressure $p$ belongs to $L^1$, then
$u=0$ and $p=0$; the original pressure is constant.
\end{lemma}
\begin{proof}
Lemma~\ref{lem:normalization} and \eqref{eq:Qq} give
$u\in L^6\cap L^\infty$ and $Q\in L^3$. Hence
\begin{equation*}
    |\nabla u|^2\in L^1,\qquad |u|^2\in L^3,
    \qquad Qu\in L^2.
\end{equation*}
For $\phi\in C_c^\infty(\R^3)$, H\"older's inequality yields
\begin{equation*}
    \begin{aligned}
        &\int|u|^2|\Delta\phi|\dd
           \le\norm{u}{L^6}^2\norm{\Delta\phi}{L^{3/2}},\\
        &\int|Q|\,|u|\,|\nabla\phi|\dd
           \le\norm{Q}{L^3}\norm{u}{L^6}\norm{\nabla\phi}{L^2}.
    \end{aligned}
\end{equation*}
Together with $\int|\nabla u|^2|\phi|\dd
\le\norm{\phi}{L^\infty}\norm{\nabla u}{L^2}^2$, these
bounds justify every term in \eqref{eq:local-energy}.
The field $u\phi$ is smooth and compactly supported. In particular,
testing the momentum equation with $u\phi$ and integrating by parts gives
\begin{equation*}
    \begin{aligned}
        &\begin{aligned}[t]
            \int(-\Delta u)\cdot u\phi\dd
               &=\int|\nabla u|^2\phi\dd
                  +\frac12\int\nabla(|u|^2)\cdot\nabla\phi\dd\\
               &=\int|\nabla u|^2\phi\dd
                  -\frac12\int|u|^2\Delta\phi\dd,
        \end{aligned}\\[3pt]
        &\begin{aligned}[t]
            \int(u\cdot\nabla)u\cdot u\phi\dd
               &=\frac12\int u\cdot\nabla(|u|^2)\phi\dd
                =-\frac12\int|u|^2u\cdot\nabla\phi\dd,
        \end{aligned}\\[3pt]
        &\begin{aligned}[t]
            \int\nabla p\cdot u\phi\dd
               &=-\int p\,\Div(u\phi)\dd
                =-\int p\,u\cdot\nabla\phi\dd.
        \end{aligned}
    \end{aligned}
\end{equation*}
Here the last two identities use $\Div u=0$.
The viscous cross term is also controlled explicitly by
\begin{equation*}
    \int|\nabla u|\,|u|\,|\nabla\phi|\dd
       \le\norm{\nabla u}{L^2}\norm{u}{L^6}
                                \norm{\nabla\phi}{L^3}<\infty.
\end{equation*}
All other integrations by parts involve smooth functions on the
compact support of $\phi$. Adding the three identities proves
\eqref{eq:local-energy}. Adding a constant to $p$ does not change
this identity, since $\int u\cdot\nabla\phi\dd=0$.

Now assume that the normalized pressure belongs to $L^1$.
By \eqref{eq:sign},
\begin{equation*}
    \norm{u}{L^2}^2\le2\norm{p}{L^1},\qquad
    \norm{Q}{L^1}\le\norm{p}{L^1},\qquad
    \norm{Qu}{L^1}\le\norm{u}{L^\infty}\norm{p}{L^1}.
\end{equation*}
Thus both $|u|^2$ and $Qu$ are globally integrable.
Choose $\theta\in C_c^\infty(B_2)$ with $0\le\theta\le1$
and $\theta=1$ on $B_1$, and put $\phi_R(x)=\theta(x/R)$.
Then $\phi_R\to1$ pointwise,
$|\nabla\phi_R|\le CR^{-1}$ and
$|\Delta\phi_R|\le CR^{-2}$. Applying \eqref{eq:local-energy},
we obtain
\begin{align*}
    0\le\int|\nabla u|^2\phi_R\dd
       &\le CR^{-2}\norm{u}{L^2}^2
          +CR^{-1}\norm{Qu}{L^1}\\
       &\le C\bigl(R^{-2}+R^{-1}\norm{u}{L^\infty}\bigr)
                                            \norm{p}{L^1}
          \longrightarrow0.
\end{align*}
Since $|\nabla u|^2\in L^1$ and $0\le\phi_R\le1$, dominated
convergence gives $\int|\nabla u|^2\dd=0$. Smoothness then implies
that $u$ is constant, and $u\in L^2$ gives $u=0$.
The momentum equation yields $\nabla p=0$, and $p\in L^1$
forces the normalized pressure to be zero.
\end{proof}
\subsection{Cut off functions}
We fix one family of cylindrical cutoffs for both component arguments.
Choose $0\le\chi\le1$ in $C_c^\infty(\R^2)$, radial, equal to one
on $|x_h|\le1$ and zero on $|x_h|\ge2$. Choose $0\le\eta\le1$
in $C_c^\infty(\R)$, equal to one on $[-1,1]$ and zero outside
$[-2,2]$. Set
\begin{equation}\label{eq:cylinder-cutoff}
    \Phi_{R,H}(x)=\chi(x_h/R)\eta(z/H),
    \qquad E_{R,H}=\{\lvert x_h\rvert<2R,\ |z|<2H\}.
\end{equation}
Then $|E_{R,H}|\le CR^2H$, and
\begin{align}\label{eq:cutoff-derivatives}
    |\nabla_h\Phi_{R,H}|
       &\le CR^{-1}\ind_{\{R<\lvert x_h\rvert<2R,\ |z|<2H\}},\notag\\
    |\partial_z\Phi_{R,H}|
       &\le CH^{-1}\ind_{\{\lvert x_h\rvert<2R,\ H<|z|<2H\}},\\
    |\Delta\Phi_{R,H}|
       &\le C(R^{-2}+H^{-2})\ind_{E_{R,H}}.\notag
\end{align}
In particular, $\Phi_{R,H}\to1$ pointwise when $R,H\to\infty$.
The difference between the two proofs will be the pressure estimate
and the choice of $H$ relative to $R$. In both proofs, a horizontal
cutoff derivative pairs with $u_h$, while a vertical derivative pairs
with $w$. Keeping $R$ and $H$ independent until the errors have been
estimated makes this distinction explicit. A ball cutoff would not
allow the two fluxes to be adjusted separately.

\section{The proof of Theorem \ref{thm:two}}\label{sec:two}
In this section, we prove Theorem~\ref{thm:two}. The main estimate controls
$\nabla p$ in $L^{4/3,2}$ and $p/\lvert x_h\rvert$ in the same space via the condition $\uh\in L^{4,\infty}$. The latter
quantity is adapted to the horizontal derivative of a cylindrical
cutoff. We first establish the endpoint tools, then derive the pressure
bounds, and finally close the energy identity.

\subsection{Cylindrical Hardy estimates}
The next three lemmas concern function-space properties only and
do not use the Navier--Stokes equations. In each lemma, $f$ denotes
an arbitrary function satisfying the stated assumptions; its meaning
is local to that lemma.

\begin{lemma}\label{lem:hardy}
Let $1<r<2$ and $1\le s<\infty$. If $\nabla_h f\in L^{r,s}$, there exists a constant $C_{r,s}$ such that
\begin{equation}\label{eq:hardy}
    \norm{f/\lvert x_h\rvert}{L^{r,s}}
       \le C_{r,s}\norm{\nabla_h f}{L^{r,s}}.
\end{equation}
\end{lemma}
\begin{proof}
By a density argument, it suffices to consider
$f\in C_c^\infty(\R^3\setminus\{x_h=0\})$.
Indeed, for a fixed smooth compactly supported function, a horizontal
cutoff that vanishes for $|x_h|\le\delta$ and equals one for
$|x_h|\ge2\delta$ produces an additional gradient term of size
$O(\delta^{2/r-1})$ in $L^{r,s}$. This tends to zero because $r<2$;
the remaining gradient term converges by
Lemma~\ref{lem:lorentz-basic}(iv).
Thus we may assume that $f$ is smooth and its support stays away
from the axis $x_h=0$.
For $|x_h|>0$, integration along a horizontal ray gives
\begin{equation*}
    f(x_h,z)=-\int_1^\infty x_h\cdot\nabla_hf(tx_h,z)\,\mathrm{d}t,
    \qquad
    \frac{|f(x_h,z)|}{\lvert x_h\rvert}
       \le\int_1^\infty|\nabla_hf(tx_h,z)|\,\mathrm{d}t.
\end{equation*}
For $t>0$, dilation of the two horizontal coordinates gives
\begin{equation*}
    \norm{(\nabla_h f)(t x_h,z)}{L^{r,s}}
       =t^{-2/r}\norm{\nabla_h f}{L^{r,s}}.
\end{equation*}
Using the equivalent Banach norm from Lemma~\ref{lem:lorentz-basic}(i),
Minkowski's inequality yields
\begin{equation*}
    \norm{f/\lvert x_h\rvert}{L^{r,s}}
       \le C_{r,s}\norm{\nabla_hf}{L^{r,s}}
                         \int_1^\infty t^{-2/r}\,\mathrm{d}t.
\end{equation*}
The integral is finite because $r<2$. The Fatou property then
removes the approximation and gives the stated extension.
The values on the axis $x_h=0$ do not matter.
\end{proof}

\begin{lemma}\label{lem:log}
Let $f\in L^{4,\infty}\cap L^\infty$ with
$\norm{f}{L^{4,\infty}}>0$. For every measurable set $E$ of finite measure, it holds that
\begin{equation}\label{eq:log}
    \norm{f\ind_E}{L^{4,2}}^2
       \le \norm{f}{L^{4,\infty}}^2
          \left[2+\logp\left(
             \frac{|E|\norm{f}{L^\infty}^4}
                  {\norm{f}{L^{4,\infty}}^4}\right)\right],
\end{equation}
where $\logp t=\max\{0,\log t\}$ for $t>0$ and $\logp0=0$.
\end{lemma}
\begin{proof}
If $|E|=0$, the inequality is trivial.
For $0<t<|E|$, the restricted rearrangement is bounded by
$\min\{\norm{f}{L^\infty},\norm{f}{L^{4,\infty}}t^{-1/4}\}$,
and it is zero for $t>|E|$. Therefore
\begin{equation*}
    \norm{f\ind_E}{L^{4,2}}^2
       \le\int_0^{|E|}
          \min\{\norm{f}{L^\infty}^2t^{1/2},
                \norm{f}{L^{4,\infty}}^2\}\,\frac{\mathrm{d}t}{t}.
\end{equation*}
Split the integral at
$t_0=(\norm{f}{L^{4,\infty}}/\norm{f}{L^\infty})^4$.
The part below $t_0$ is at most $2\norm{f}{L^{4,\infty}}^2$;
the remaining part is
$\norm{f}{L^{4,\infty}}^2\logp(|E|/t_0)$.
This proves \eqref{eq:log}.
\end{proof}

\begin{lemma}\label{lem:annuli}
Let $f\in L^{4/3,2}$, and define
\begin{equation*}
    \mathcal A_n=\{2^n<\lvert x_h\rvert<2^{n+1}\},\qquad n\ge1.
\end{equation*}
Then we have
\begin{equation}\label{eq:square-sum}
    \sum_{n=1}^\infty \norm{f\ind_{\mathcal A_n}}{L^{4/3,2}}^2\le\norm{f}{L^{4/3,2}}^2.
\end{equation}
Consequently, there is a strictly increasing sequence $n_j\to\infty$
such that
\begin{equation}\label{eq:subsequence}
    \sqrt{n_j+1}\,\norm{f\ind_{\mathcal A_{n_j}}}{L^{4/3,2}}\longrightarrow0.
\end{equation}
\end{lemma}
\begin{proof}
By the definition of $L^{4/3,2}$, we see that
\begin{equation}\label{eq:layercake}
    \norm{f}{L^{4/3,2}}^2
       =4/3\int_0^\infty\lambda d_f(\lambda)^{3/2}\,\mathrm{d}\lambda.
\end{equation}
Since the annuli are disjoint and $3/2>1$,
\begin{equation*}
    \sum_n d_{f\ind_{\mathcal A_n}}(\lambda)^{3/2}
       \le\left(\sum_n d_{f\ind_{\mathcal A_n}}(\lambda)\right)^{3/2}
       \le d_f(\lambda)^{3/2}.
\end{equation*}
Combining this inequality with \eqref{eq:layercake} proves
\eqref{eq:square-sum}. If \eqref{eq:subsequence} failed for every
subsequence, then $(n+1)\norm{f\ind_{\mathcal A_n}}{L^{4/3,2}}^2$ would have a positive lower bound for
all sufficiently large $n$. This would contradict
\eqref{eq:square-sum} and the divergence of $\sum_n(n+1)^{-1}$.
Equivalently, square-summability implies
\begin{equation*}
    \liminf_{n\to\infty}(n+1)\norm{f\ind_{\mathcal A_n}}{L^{4/3,2}}^2=0.
\end{equation*}
This provides a subsequential rate that does not follow from
$\norm{f\ind_{\mathcal A_n}}{L^{4/3,2}}\to0$ alone. Mere convergence gives no rate capable of absorbing
a factor $\sqrt{n+1}$. The subsequence is chosen using $f$ on entire
cylindrical annuli, so its choice will not depend on the cylinder heights.
\end{proof}

\subsection{Pressure estimates from incompressibility}
\begin{proposition}\label{prop:two-pressure}
Under \eqref{eq:NS}--\eqref{eq:background}, assume
$\uh\in L^{4,\infty}$. Then
\begin{align}
    \norm{\nabla p}{L^{4/3,2}}&\le C\norm{\uh}{L^{4,\infty}}\norm{\nabla u}{L^2},\label{eq:pressure-gradient}\\
    \norm{p}{L^{12/5,2}}&\le C\norm{\uh}{L^{4,\infty}}\norm{\nabla u}{L^2},\label{eq:pressure-sobolev}\\
    \norm{p/\lvert x_h\rvert}{L^{4/3,2}}&\le C\norm{\uh}{L^{4,\infty}}\norm{\nabla u}{L^2}.\label{eq:pressure-hardy}
\end{align}
\end{proposition}
\begin{proof}
We first use incompressibility in the pressure equation, and then justify
the two homogeneous estimates.
\step{Step 1. Cancellation of the purely vertical term.}
Incompressibility gives
\begin{equation*}
    \partial_z(w^2)
       =2w\partial_z w
       =-2w\Div_h\uh
       =-2\Div_h(\uh w)+2\uh\cdot\nabla_hw.
\end{equation*}
Using this identity in \eqref{eq:poisson}, we obtain
\begin{align}\label{eq:pressure-cancellation}
    -\Delta p
       &=\sum_{\alpha,\beta=1}^2
           \partial_\alpha\partial_\beta(u_\alpha u_\beta)
          +2\sum_{\alpha=1}^2\partial_\alpha\partial_z(u_\alpha w)
          +\partial_z^2(w^2)\notag\\
       &=\sum_{\alpha,\beta=1}^2
           \partial_\alpha\partial_\beta(u_\alpha u_\beta)
          +2\partial_z(\uh\cdot\nabla_hw).
\end{align}
Consequently,
\begin{equation}\label{eq:pressure-representation}
    p=\sum_{\alpha,\beta=1}^2\RR_\alpha\RR_\beta(u_\alpha u_\beta)
         +2\partial_z(-\Delta)^{-1}(\uh\cdot\nabla_hw).
\end{equation}
To justify the representation, the first term belongs to $L^3$.
Also, $\uh\cdot\nabla_hw\in L^{3/2}$ by $u\in L^6$ and
$\nabla u\in L^2$, so its order-one potential belongs to $L^3$.
The difference between the two sides of \eqref{eq:pressure-representation}
is therefore an entire harmonic $L^3$ function, and hence zero.

\step{Step 2. The pressure gradient.}
Differentiate \eqref{eq:pressure-representation} in distributions:
\begin{equation}\label{eq:gradient-representation}
\begin{aligned}
    \partial_kp
       ={}&\sum_{\alpha,\beta=1}^2\RR_\alpha\RR_\beta
          \bigl((\partial_ku_\alpha)u_\beta
                           +u_\alpha\partial_ku_\beta\bigr)\\
        &+2\RR_k\RR_3
                                  (\uh\cdot\nabla_hw).
\end{aligned}
\end{equation}
Every operator on the right is of order zero, and every product is
bounded by a constant times $|\uh||\nabla u|$. By \eqref{eq:product},
\begin{equation*}
    \norm{|\uh||\nabla u|}{L^{4/3,2}}
       \le C\norm{\uh}{L^{4,\infty}}\norm{\nabla u}{L^2}.
\end{equation*}
Here $L^{2,2}=L^2$ for the normalization in
\eqref{eq:lorentz-definition}. Boundedness of the singular integrals
on $L^{4/3,2}$ proves \eqref{eq:pressure-gradient}. The finite second
Lorentz index is supplied by the Dirichlet integral, and will be
retained in the weighted pressure estimate and the annular square sum.

\step{Step 3. The unweighted pressure.}
Define
\begin{equation*}
    \widetilde p=-\sum_{j=1}^3
                   \partial_j(-\Delta)^{-1}(\partial_jp).
\end{equation*}
The potential estimate \eqref{eq:potential} gives
\begin{equation*}
    \norm{\widetilde p}{L^{12/5,2}}
       \le C\norm{\nabla p}{L^{4/3,2}},
    \qquad \frac{1}{12/5}=\frac{1}{4/3}-\frac13.
\end{equation*}
The difference $h=p-\widetilde p$ is harmonic in distributions,
and Weyl's lemma gives it a smooth harmonic representative. For every
fixed $x$ and every $R>0$, the mean-value property gives
\begin{align*}
    |h(x)|
       &\le\frac{1}{|B_R|}\int_{B_R(x)}
                          (|p|+|\widetilde p|)\dd\\
       &\le C\left(R^{-1}\norm{p}{L^3}
                       +R^{-5/4}\norm{\widetilde p}{L^{12/5,2}}\right).
\end{align*}
For the first term, H\"older's inequality contributes
$|B_R|^{-1/3}=CR^{-1}$. For the second, use
$L^{12/5,2}\subset L^{12/5,\infty}$ and
\eqref{eq:finite-volume}; the factor is
$|B_R|^{-5/12}=CR^{-5/4}$.
Letting $R\to\infty$ proves $p=\widetilde p$ and
\eqref{eq:pressure-sobolev}. This also fixes the representative in the
homogeneous estimate.

\step{Step 4. The cylindrical weighted pressure.}
Let $\theta_L$ be a standard full-space cutoff, equal to one on $B_L$
and supported in $B_{2L}$, with $|\nabla\theta_L|\le C/L$.
Then
\begin{equation*}
    \nabla_h(p\theta_L)-\nabla_hp
       =(\theta_L-1)\nabla_hp+p\nabla_h\theta_L.
\end{equation*}
The first term tends to zero in $L^{4/3,2}$ by absolute continuity.
For the second, \eqref{eq:product} gives
\begin{align*}
    \norm{p\nabla_h\theta_L}{L^{4/3,2}}
       &\le C\norm{p\ind_{\{L<|x|<2L\}}}{L^{12/5,2}}
                       \norm{\nabla_h\theta_L}{L^{3,\infty}}\\
       &\le C\norm{p\ind_{\{L<|x|<2L\}}}{L^{12/5,2}}
          \longrightarrow0.
\end{align*}
The cutoff-gradient norm is uniformly bounded, and the second Lorentz
index of $p$ is finite. Since $p$ is smooth, $p\theta_L$ is an
admissible compactly supported function in Lemma~\ref{lem:hardy}.
Passing to the limit proves \eqref{eq:pressure-hardy}.
\end{proof}

\subsection{The energy flux on elongated cylinders}
\begin{proof}[Proof of Theorem~\ref{thm:two}]
By \eqref{eq:two-condition}, $\uh\in L^{4,\infty}$.
Proposition~\ref{prop:two-pressure} gives
$p/\lvert x_h\rvert\in L^{4/3,2}$ and $p\in L^{12/5,2}$.
If $\norm{\uh}{L^{4,\infty}}=0$, \eqref{eq:pressure-gradient}
gives $\nabla p=0$. Since $p\in L^3$, we have $p=0$, and
\eqref{eq:sign} implies $u=0$.
We therefore assume $\norm{\uh}{L^{4,\infty}}>0$.

\step{Step 1. The three terms in the energy identity.}
Use $\phi=\Phi_{R,H}$ from \eqref{eq:cylinder-cutoff} in
\eqref{eq:local-energy}. Write
\begin{equation}\label{eq:two-energy}
    \int|\nabla u|^2\Phi_{R,H}\dd
       =J_1+J_2+J_3,
\end{equation}
where
\begin{align*}
    J_1&=\frac12\int|u|^2\Delta\Phi_{R,H}\dd,\\
    J_2&=\int Q\uh\cdot\nabla_h\Phi_{R,H}\dd,\qquad
    J_3=\int Qw\partial_z\Phi_{R,H}\dd.
\end{align*}
We estimate the diffusion and vertical terms using the two pressure
norms, and then select radii for the horizontal flux.

\step{Step 2. The diffusion term.}
On $E_{R,H}$, \eqref{eq:sign} gives $|u|^2\le2q$, while
$q\le2R(q/\lvert x_h\rvert)$. By \eqref{eq:finite-volume},
\begin{equation*}
    \int_{E_{R,H}}q\dd
       \le CR\norm{p/\lvert x_h\rvert}{L^{4/3,2}}|E_{R,H}|^{1/4}
       \le C\norm{p/\lvert x_h\rvert}{L^{4/3,2}}R^{3/2}H^{1/4}.
\end{equation*}
Using \eqref{eq:cutoff-derivatives}, we obtain
\begin{equation}\label{eq:two-diffusion}
    |J_1|
       \le C\norm{p/\lvert x_h\rvert}{L^{4/3,2}}\left(R^{-1/2}H^{1/4}+R^{3/2}H^{-7/4}\right).
\end{equation}

\step{Step 3. The vertical flux.}
The sign in \eqref{eq:sign} gives $|Qw|\le\sqrt{2}\,q^{3/2}$.
The finite-volume estimate with $r=12/5$ and $a=3/2$ yields
\begin{equation*}
    \int_{E_{R,H}}q^{3/2}\dd
       \le C\norm{p}{L^{12/5,2}}^{3/2}|E_{R,H}|^{3/8}.
\end{equation*}
Therefore
\begin{equation}\label{eq:two-vertical}
    |J_3|\le C\norm{p}{L^{12/5,2}}^{3/2}R^{3/4}H^{-5/8}.
\end{equation}
No separate integrability assumption on $w$ is used here.

\step{Step 4. The horizontal flux at the weak endpoint.}
Let $R=2^n$ and $H=R^{3/2}$, and let $\mathcal A_n$ be the
annuli in Lemma~\ref{lem:annuli}. On the support of
$\nabla_h\Phi_{R,H}$, we have $R<\lvert x_h\rvert<2R$.
Thus \eqref{eq:sign} and Lemma~\ref{lem:lorentz-holder} imply
\begin{align}\label{eq:two-horizontal}
    |J_2|
       &\le C\int_{\mathcal A_n\cap E_{R,H}}
                       \frac{q}{\lvert x_h\rvert}|\uh|\dd\notag\\
       &\le C\norm{(p/\lvert x_h\rvert)\ind_{\mathcal A_n}}{L^{4/3,2}}
           \norm{\uh\ind_{\mathcal A_n\cap E_{R,H}}}{L^{4,2}}.
\end{align}
The second factor is estimated locally. Since
$|E_{R,H}|\le C2^{7n/2}$, Lemma~\ref{lem:log} gives
\begin{equation*}
\begin{aligned}
    \norm{\uh\ind_{\mathcal A_n\cap E_{R,H}}}{L^{4,2}}
       \le{}& C\norm{\uh}{L^{4,\infty}}
          \left[1+\logp\left(
             \frac{\norm{\uh}{L^\infty}^4}
                  {\norm{\uh}{L^{4,\infty}}^4}\right)\right]^{1/2}
          \sqrt{n+1}.
\end{aligned}
\end{equation*}
Consequently,
\begin{equation}\label{eq:two-horizontal-final}
\begin{aligned}
    |J_2|\le{}& C\norm{\uh}{L^{4,\infty}}
          \left[1+\logp\left(
             \frac{\norm{\uh}{L^\infty}^4}
                  {\norm{\uh}{L^{4,\infty}}^4}\right)\right]^{1/2}\\
       &\quad\times\sqrt{n+1}\,
          \norm{(p/\lvert x_h\rvert)\ind_{\mathcal A_n}}{L^{4/3,2}}.
\end{aligned}
\end{equation}
The first line is independent of $n$. Lemma~\ref{lem:annuli},
applied to $f=p/\lvert x_h\rvert$, supplies an increasing sequence $n_j$
along which the second line tends to zero. The logarithm in
Lemma~\ref{lem:log} grows at most linearly in $n$ because the cylinder
volume grows like a fixed power of $2^n$. Its square root therefore
matches the factor removed by the annular square sum. No estimate on
the rate of decay of a weak-$L^4$ tail is needed.

\step{Step 5. A common sequence for all fluxes.}
Set $R_j=2^{n_j}$ and $H_j=R_j^{3/2}$. The preceding bounds give
\begin{align*}
    |J_1|&\le C\norm{p/\lvert x_h\rvert}{L^{4/3,2}}\left(R_j^{-1/8}+R_j^{-9/8}\right)\longrightarrow0,\\
    |J_3|&\le C\norm{p}{L^{12/5,2}}^{3/2}R_j^{-3/16}\longrightarrow0.
\end{align*}
Also, $J_2\to0$ along this sequence by
\eqref{eq:two-horizontal-final} and Lemma~\ref{lem:annuli}.
Since $0\le\Phi_{R_j,H_j}\le1$, the cutoffs tend pointwise to one,
and $|\nabla u|^2\in L^1$, dominated convergence in
\eqref{eq:two-energy} gives $\int|\nabla u|^2\dd=0$.
The decay condition implies $u=0$, and \eqref{eq:NS} then gives
$\nabla p=0$.
\end{proof}

\begin{remark}\label{rem:two-height}
The height exponent $3/2$ is convenient but not essential. For
$H=R^a$, the three powers in \eqref{eq:two-diffusion} and
\eqref{eq:two-vertical} are
\begin{equation*}
    -\frac12+\frac a4,\qquad
    \frac32-\frac{7a}{4},\qquad
    \frac34-\frac{5a}{8}.
\end{equation*}
They are all negative when $6/5<a<2$. The cylinder volume is then
bounded by $CR^{2+a}$, so the logarithmic estimate remains of order
$\sqrt{n+1}$. In contrast, $H=R$ leaves the positive power $R^{1/8}$
in the vertical flux bound. The annular selection only needs to hold
along a subsequence; convergence of the weak-$L^4$ tail is not assumed.
\end{remark}

\section{The proof of Theorem~\ref{thm:one}}\label{sec:one}
In this section we prove Theorem~\ref{thm:one}. The pressure-gradient
estimate of Section~\ref{sec:two} is invalid. Instead, we test the
third momentum equation with a bounded approximation to the coordinate
$z$ and show that the total normalized pressure is integrable.

Firstly, we prove a standard interpolation estimate.
\begin{lemma}\label{lem:interpolation}
If $f\in L^{3/2,\infty}\cap L^6$, then
\begin{equation}\label{eq:interpolation}
    \norm{f}{L^2}
       \le C\norm{f}{L^{3/2,\infty}}^{2/3}
              \norm{f}{L^6}^{1/3}.
\end{equation}
\end{lemma}
\begin{proof}
The assertion is immediate if $f=0$ almost everywhere. Otherwise,
\begin{equation*}
    d_f(\lambda)
       \le\min\{\norm{f}{L^{3/2,\infty}}^{3/2}\lambda^{-3/2},
                 \norm{f}{L^6}^6\lambda^{-6}\}.
\end{equation*}
Split $2\int_0^\infty\lambda d_f(\lambda)\,\mathrm{d}\lambda$ at
\begin{equation*}
    \lambda_0=\norm{f}{L^6}^{4/3}\norm{f}{L^{3/2,\infty}}^{-1/3}.
\end{equation*}
Both resulting terms are bounded by
$C\norm{f}{L^{3/2,\infty}}^{4/3}\norm{f}{L^6}^{2/3}$,
which proves \eqref{eq:interpolation}.
\end{proof}
We introduce some new test functions.
For $\ell\ge1$, define
\begin{equation}\label{eq:psi}
    \psi_\ell(x)=\frac{\ell z}{\sqrt{\ell^2+|x|^2}},
    \qquad a_\ell=\partial_z\psi_\ell,
    \qquad b_\ell=\nabla_h\psi_\ell.
\end{equation}
\Needspace{14\baselineskip}
\begin{lemma}\label{lem:weight}
The functions in \eqref{eq:psi} satisfy
\begin{equation}\label{eq:weight-derivatives}
    a_\ell=\frac{\ell(\ell^2+\lvert x_h\rvert^2)}{(\ell^2+|x|^2)^{3/2}},
    \qquad b_\ell=-\frac{\ell z x_h}{(\ell^2+|x|^2)^{3/2}},
\end{equation}
with
\begin{equation}\label{eq:weight-bounds}
    0<a_\ell\le1,\qquad \frac{|b_\ell|^2}{a_\ell}\le1,
    \qquad a_\ell(x)\longrightarrow1\quad(\ell\to\infty).
\end{equation}
Moreover,
\begin{equation}\label{eq:weight-uniform}
    |\psi_\ell|\le\ell,\qquad
    |\nabla\psi_\ell|\le\frac{C\ell}{\sqrt{\ell^2+|x|^2}},
    \qquad \norm{\Delta\psi_\ell}{L^{3,1}}\le C.
\end{equation}
For the solution under consideration,
\begin{equation}\label{eq:absorption}
    w\,u\cdot\nabla\psi_\ell
       \le\frac12q a_\ell+2w^2.
\end{equation}
\end{lemma}
\begin{proof}
Direct differentiation gives \eqref{eq:weight-derivatives}. In particular,
\begin{equation*}
    \frac{|b_\ell|^2}{a_\ell}
       =\frac{\ell z^2\lvert x_h\rvert^2}
              {(\ell^2+|x|^2)^{3/2}(\ell^2+\lvert x_h\rvert^2)}
       \le\frac{\ell z^2}{(\ell^2+|x|^2)^{3/2}}\le1.
\end{equation*}
The remaining pointwise statements follow from the same formulas.
A further differentiation gives
\begin{equation}\label{eq:laplace-psi}
    \Delta\psi_\ell(x)
       =-\frac{\ell z(5\ell^2+2|x|^2)}
                      {(\ell^2+|x|^2)^{5/2}}.
\end{equation}
Hence $|\Delta\psi_1(x)|\le C(1+|x|)^{-2}$, whose rearrangement is
bounded by $C\min\{1,t^{-2/3}\}$. This proves
$\Delta\psi_1\in L^{3,1}$. Since
$\Delta\psi_\ell(x)=\ell^{-1}(\Delta\psi_1)(x/\ell)$,
the $L^{3,1}$ norm is independent of the scale $\ell$.

Finally, \eqref{eq:sign} implies $|\uh|\le\sqrt{2q}$. Young's
inequality gives
\begin{equation*}
    |w\uh\cdot b_\ell|
       \le |w|\sqrt{2q}\,|b_\ell|
       \le\frac12q a_\ell+w^2\frac{|b_\ell|^2}{a_\ell}.
\end{equation*}
Adding $w^2a_\ell$ and using \eqref{eq:weight-bounds} proves
\eqref{eq:absorption}.
\end{proof}

\subsection{Recovery of integrable pressure}
The scalar weak form of the third momentum equation is
\begin{equation}\label{eq:scalar-weak}
    \int q\,\partial_z\varphi\dd
       =\int w\,\Delta\varphi\dd
          +\int wu\cdot\nabla\varphi\dd,
    \qquad \varphi\in C_c^\infty(\R^3).
\end{equation}
The positivity of $a_\ell$ is useful here, but the test $\psi_\ell$
is not compactly supported. We justify it before drawing any conclusion
about $\int q\dd$.

\begin{proposition}\label{prop:one-pressure}
Under \eqref{eq:NS}--\eqref{eq:background}, assume
$w\in L^{3/2,\infty}$. Then
\begin{equation}\label{eq:one-pressure}
    \int q\dd=\norm{p}{L^1}
       \le C\left(\norm{w}{L^{3/2,\infty}}+\norm{w}{L^2}^2\right).
\end{equation}
In particular, $p\in L^1$ and $u\in L^2$.
\end{proposition}
\begin{proof}
Lemma~\ref{lem:interpolation} gives $w\in L^2$.

\step{Step 1. A compactly supported identity.}
For $R,H\ge\ell$, use
$\varphi=\psi_\ell\Phi_{R,H}$ in \eqref{eq:scalar-weak}, and put
\begin{equation*}
    I_{\ell,R,H}=\int q a_\ell\Phi_{R,H}\dd.
\end{equation*}
Expansion of the derivatives yields the exact identity
\begin{equation}\label{eq:one-expanded}
\begin{aligned}
    I_{\ell,R,H}
       ={}&\int w\Phi_{R,H}\Delta\psi_\ell\dd
          +\int wu\cdot\nabla\psi_\ell\,\Phi_{R,H}\dd\\
        &+\int w\bigl(2\nabla\psi_\ell\cdot\nabla\Phi_{R,H}
                           +\psi_\ell\Delta\Phi_{R,H}\bigr)\dd\\
        &+\int wu\cdot\nabla\Phi_{R,H}\,\psi_\ell\dd
          -\int q\psi_\ell\partial_z\Phi_{R,H}\dd.
\end{aligned}
\end{equation}
By Lemma~\ref{lem:lorentz-holder} and \eqref{eq:weight-uniform}, the absolute
value of the first integral is at most $C\norm{w}{L^{3/2,\infty}}$. The second integral is
at most $I_{\ell,R,H}/2+2\norm{w}{L^2}^2$ by \eqref{eq:absorption}.
Let $E_1$, $E_2$, and $E_3$
denote the integrals of the absolute values of the last three
integrands in \eqref{eq:one-expanded}. Absorption gives
\begin{equation}\label{eq:one-master}
    I_{\ell,R,H}
       \le C(\norm{w}{L^{3/2,\infty}}+\norm{w}{L^2}^2)
          +C(E_1+E_2+E_3).
\end{equation}

\step{Step 2. Viscous cutoff errors.}
The horizontal and vertical derivatives of $\Phi_{R,H}$ are supported,
respectively, in
\begin{equation*}
    \mathcal A_h=\{R\le\lvert x_h\rvert\le2R,\ |z|\le2H\},
    \qquad \mathcal A_z=\{\lvert x_h\rvert\le2R,\ H\le|z|\le2H\}.
\end{equation*}
Both sets have measure at most $CHR^2$.
On $\mathcal A_h$, $|\nabla\psi_\ell|\le C\ell/R$, and on $\mathcal A_z$,
$|\nabla\psi_\ell|\le C\ell/H$. Consequently,
\begin{equation*}
    E_1
       \le C\ell R^{-2}\int_{\mathcal A_h}|w|\dd
          +C\ell H^{-2}\int_{\mathcal A_z}|w|\dd.
\end{equation*}
Using \eqref{eq:finite-volume} with $r=3/2$ and $a=1$, we obtain
\begin{equation}\label{eq:one-viscous}
    E_1
       \le C\ell \norm{w}{L^{3/2,\infty}}\left(H^{1/3}R^{-4/3}
                            +H^{-5/3}R^{2/3}\right).
\end{equation}

\step{Step 3. Convective cutoff errors.}
The vertical contribution contains $w^2$, so
\begin{equation*}
    \int|w|^2|\partial_z\Phi_{R,H}\,\psi_\ell|\dd\le C\ell H^{-1}\norm{w}{L^2}^2.
\end{equation*}
For the horizontal contribution, H\"older's inequality and
\eqref{eq:finite-volume} with $a=6/5$ give
\begin{align*}
    \int|w\uh\cdot\nabla_h\Phi_{R,H}\,\psi_\ell|\dd
       &\le C\ell R^{-1}\norm{w}{L^{6/5}(\mathcal A_h)}
                              \norm{u}{L^6(\mathcal A_h)}\\
       &\le C\ell \norm{w}{L^{3/2,\infty}} H^{1/6}R^{-2/3}
                              \norm{u}{L^6(\{\lvert x_h\rvert>R\})}.
\end{align*}
Thus
\begin{equation}\label{eq:one-convective}
    E_2
       \le C\ell H^{-1}\norm{w}{L^2}^2
          +C\ell \norm{w}{L^{3/2,\infty}} H^{1/6}R^{-2/3}
                              \norm{u}{L^6(\{\lvert x_h\rvert>R\})}.
\end{equation}

\step{Step 4. The pressure cutoff error.}
Since
$p\in L^3$, it follows that
\begin{align}\label{eq:one-pressure-error}
    E_3
       &\le C\ell H^{-1}\int_{\mathcal A_z}|p|\dd\notag\\
       &\le C\ell H^{-1/3}R^{4/3}
                              \norm{p}{L^3(\{|z|>H\})}.
\end{align}

\step{Step 5. Three ordered limits.}
The powers in the preceding bounds explain the choice of cylinder.
Writing $H=R^a$, the horizontal viscous factor is $R^{(a-4)/3}$,
the horizontal convective factor is $R^{(a-4)/6}$, and the pressure
factor is $R^{(4-a)/3}$. Since we have no prescribed rate for the
strong $L^6$ and $L^3$ tails, these estimates lead to the balance
$a=4$. At this value the transverse viscous error has no vanishing
tail attached to it. An additional fixed parameter will remove it.

Fix $\ell\ge1$ and $0<\delta\le1$, and set
\begin{equation}\label{eq:one-height}
    H=\delta R^4.
\end{equation}
For all sufficiently large $R$, the restrictions $R,H\ge\ell$ hold.
Since $u\in L^6$ and $p\in L^3$,
\begin{equation*}
    \norm{u}{L^6(\{\lvert x_h\rvert>R\})}\longrightarrow0\quad(R\to\infty),
    \qquad
    \norm{p}{L^3(\{|z|>H\})}\longrightarrow0\quad(H\to\infty).
\end{equation*}
Substituting \eqref{eq:one-height} into the error estimates gives
\begin{align}
    E_1
       &\le C\ell \norm{w}{L^{3/2,\infty}}\left(\delta^{1/3}
                             +\delta^{-5/3}R^{-6}\right),
                                      \label{eq:one-scaled-viscous}\\
    E_2
       &\le C\ell \norm{w}{L^2}^2\delta^{-1}R^{-4}
          +C\ell \norm{w}{L^{3/2,\infty}}\delta^{1/6}\norm{u}{L^6(\{\lvert x_h\rvert>R\})},
                                      \label{eq:one-scaled-convective}\\
    E_3
       &\le C\ell\delta^{-1/3}\norm{p}{L^3(\{|z|>\delta R^4\})}.
                                      \label{eq:one-scaled-pressure}
\end{align}
First let $R\to\infty$ with $\ell$ and $\delta$ fixed.
Since $q a_\ell\Phi_{R,\delta R^4}\ge0$ and the cutoff tends
pointwise to one, Fatou's lemma and \eqref{eq:one-master} imply
\begin{equation*}
    \int q a_\ell\dd
       \le C(\norm{w}{L^{3/2,\infty}}+\norm{w}{L^2}^2)+C\ell \norm{w}{L^{3/2,\infty}}\delta^{1/3}.
\end{equation*}
Next let $\delta\downarrow0$, keeping $\ell$ fixed. This gives
\begin{equation}\label{eq:one-weighted-pressure}
    \int q a_\ell\dd\le C(\norm{w}{L^{3/2,\infty}}+\norm{w}{L^2}^2)
    \qquad\text{for every }\ell\ge1.
\end{equation}
Finally, let $\ell\to\infty$. Since $a_\ell(x)\to1$, another
application of Fatou's lemma proves \eqref{eq:one-pressure}.
Neither application requires monotonicity of the cutoffs or of
$a_\ell$; nonnegativity and pointwise convergence suffice. The order
of the limits is essential for these estimates. With $\delta$ fixed,
the factor $\delta^{-1/3}$ does not affect the vanishing pressure
tail. Once $R$ has gone to infinity, keeping $\ell$ fixed allows the
remaining $\ell\delta^{1/3}$ term to vanish. The resulting bound
\eqref{eq:one-weighted-pressure} is uniform in $\ell$.
The $L^2$ conclusion follows from $|u|^2\le2q$.
\end{proof}

\begin{proof}[Proof of Theorem~\ref{thm:one}]
The controlled component is $w=u_3\in L^{3/2,\infty}$.
Proposition~\ref{prop:one-pressure} gives $p\in L^1$, and
Lemma~\ref{lem:energy} proves the conclusion.
\end{proof}

\begin{remark}\label{rem:one-endpoint}
We only use the  $L^6$ tail of $u$ and the $L^3$ tail of $p$
are tend to zero. The residual transverse viscous error is removed
by the separate limit $\delta\downarrow0$. In particular, the three
limits in Step~5 are ordered; no rate of decay for $\norm{u}{L^6(\{\lvert x_h\rvert>R\})}$ or $\norm{p}{L^3(\{|z|>H\})}$ is
assumed.
\end{remark}

\section{Consequences and endpoint remarks}\label{sec:consequences}
\subsection{Smaller exponents and conditions at infinity}
As with the main theorems, the following statements hold in any
fixed orthonormal coordinate system.
\begin{corollary}\label{cor:exterior}
Let $(u,p)$ be a smooth solution of \eqref{eq:NS} satisfying
\eqref{eq:background}. Then $u=0$ and $p$ is constant if, for some
$R_0\ge0$, either
\begin{equation}\label{eq:exterior-one}
    u_3\ind_{\{|x|>R_0\}}\in L^{a,\infty}(\R^3),
    \qquad 1\le a\le\frac32,
\end{equation}
or
\begin{equation}\label{eq:exterior-two}
    \uh\ind_{\{|x|>R_0\}}\in L^{a,\infty}(\R^3),
    \qquad 1\le a\le4.
\end{equation}
\end{corollary}
\begin{proof}
For a bounded measurable $f$ and $1\le a<b<\infty$,
\begin{equation}\label{eq:bounded-inclusion}
    \norm{f}{L^{b,\infty}}
       \le\norm{f}{L^\infty}^{1-a/b}
              \norm{f}{L^{a,\infty}}^{a/b}.
\end{equation}
Indeed, $d_f(\lambda)=0$ for $\lambda>\norm{f}{L^\infty}$,
and otherwise
$\lambda^b d_f(\lambda)\le
\norm{f}{L^\infty}^{b-a}\norm{f}{L^{a,\infty}}^a$.
Apply \eqref{eq:bounded-inclusion} to the exterior restriction,
with $b=3/2$ or $b=4$. The restriction to the fixed ball is bounded
and supported on a set of finite measure, and hence belongs to the
required weak space. Apply Theorem~\ref{thm:one} or
Theorem~\ref{thm:two} to the sum of the two restrictions.
\end{proof}

\begin{corollary}\label{cor:pointwise}
Under the same background assumptions, $u=0$ and $p$ is constant
if, for some $R_0>0$ and a finite constant $C_0$, either
\begin{equation}\label{eq:decay-one}
    |u_3(x)|\le C_0(1+|x|)^{-2}
    \qquad (|x|>R_0),
\end{equation}
or
\begin{equation}\label{eq:decay-two}
    |\uh(x)|\le C_0(1+|x|)^{-3/4}
    \qquad (|x|>R_0).
\end{equation}
\end{corollary}
\begin{proof}
For $b>0$, $(1+|x|)^{-3/b}\in L^{b,\infty}(\R^3)$:
its distribution function is bounded by $C\lambda^{-b}$.
Use $b=3/2$ or $b=4$ and apply Corollary~\ref{cor:exterior}.
\end{proof}

The boundedness assumption in \eqref{eq:bounded-inclusion} is important.
On the infinite-measure space $\R^3$, changing the first Lorentz index
does not by itself give a global embedding. Here boundedness follows
from the common background assumptions. The component reduction is a
separate issue: controlling one scalar projection in a smaller weak
space does not, by an embedding alone, supply a condition on a second
independent component.

\subsection{The weak endpoints and the limits of the estimates}
The weak assumptions are genuinely weaker than their strong-space
counterparts at the same exponent. For $b=3/2$ or $b=4$, the function
\begin{equation*}
    f_b(x)=(1+|x|^2)^{-3/(2b)}
\end{equation*}
belongs to $L^{b,\infty}(\R^3)$ but not to $L^b(\R^3)$.
The first assertion follows from its distribution function. For the
second, $|f_b(x)|^b\sim|x|^{-3}$ at infinity, so the radial integral
diverges logarithmically. These are function-space examples, not
nonzero stationary solutions.

The two endpoint restrictions arise from different estimates. In the
one-component proof, replacing $L^{3/2,\infty}$ by $L^{a,\infty}$
would bound the horizontal part of $E_1$ by
\begin{equation*}
    C\ell\norm{w}{L^{a,\infty}}H^{1-1/a}R^{-2/a}.
\end{equation*}
For $H=\delta R^4$, the factor involving the cutoff scales becomes
$\delta^{1-1/a}R^{4-6/a}$. At $a=3/2$, the power of $R$ vanishes
and the separate limit in $\delta$ removes the error. For $a>3/2$,
this estimate contains a positive power of $R$.

In the two-component proof, the product $|\uh||\nabla u|$ under an
$L^{a,\infty}$ assumption has first Lorentz index $r$ given by
$1/r=1/a+1/2$. Pairing the resulting weighted pressure with $\uh$
without a volume factor requires
\begin{equation*}
    \frac1r+\frac1a=1,
    \qquad\text{equivalently}\qquad
    \frac2a+\frac12=1.
\end{equation*}
This gives $a=4$. For $a>4$, the finite-volume bound contains the
positive factor $|E_{R,H}|^{1/2-2/a}$. Square-summability on dyadic
annuli does not provide a prescribed polynomial rate in $R$, so the
subsequence argument used here does not eliminate this factor.
These observations describe the scope of the two proofs; they do not
establish optimality of either component criterion.

\smallskip
\noindent\textbf{Acknowledgments} This work was partially supported by the National Key R\&D Program of China(No. 2021YFA1002100) and the National Natural Science Foundation of China (No.12571261 and No. 12471233).

\noindent\textbf{Conflict of interest}: This work does not have any conflicts of interest.

\noindent\textbf{Data Availability} No data, models, or code were generated or used during the study.

\bibliographystyle{abbrv}
\bibliography{references}

\end{document}